\documentclass[10pt]{article}
\usepackage{graphicx}
\usepackage{todonotes}

\def\R{{\rm I\! R}}

\def\C{\mbox{l\hspace{-.47em}C}}

\def \pTW*{\partial_{W^*} T}
\def \pTnW*{\partial^{(n)}_{W^*} T}

\newtheorem{theorem}{Theorem}
\newtheorem{corollary}{Corollary}
\newtheorem{definition}{Definition}
\newtheorem{lemma}{Lemma}

\title{Centers with equal period functions} 

\author{Marco Sabatini 
\footnote{Dip. di Matematica, Univ. di Trento, I-38050 Povo, (TN) - Italy.
Phone: ++39(0461)281670, Fax: ++39(0461)281624, Email: marco.sabatini@unitn.it - \ \ \ \ \ \ \ \  \ \ \ \ \ \ \ \ \ \ \ \ This paper was partially supported by the GNAMPA, {\it Gruppo Nazionale per l'Analisi Matematica, la Probabilit\`a e le loro Applicazioni}. }
}
\date{January 28th, 2015 }
\begin{document}
\maketitle
\begin{abstract}  
We give a sufficient condition for systems with symmetries to have periodic solutions with equal periods.  We show that the main result can be applied both to Hamiltonian and to non-Hamiltonian systems. We apply the main results to produce planar centers with equal period functions. 
{\bf Keywords}:  center, reversibility, symmetry, period function, boundary value problem \end{abstract}



\section{Introduction}

Let us consider a  differential system
\begin{equation} \label{sysn}
\dot z = V(z),
\end{equation}  
where $V=(V_1,\dots,V_n) \in C^1 (\Omega,\R^n)$,  $\Omega \subset \R^n$ open and connected, $z=(z_1,\dots,z_n) \in \Omega$. We denote by $\phi(t,z)$ the solution of  (\ref{sysn}) such that $\phi(0,z) = z$. Following the terminology and notation of \cite{BS}, $\phi(t,z)$ is a local flow in $\Omega$. If there exists $T > 0$ such that $\phi(T,z) = z$, then $z$ is said to be a    periodic point and its orbit $\gamma_z$ is  a cycle, support of the periodic solution $\gamma_z(t) \equiv \phi(t,z)$. The existence of periodic solutions has been for a long time one the main ones in mathematical physics, originating form the study of celestial mechanics, and subsequently extended to other classes of models in physics, biology, economics. In the simplest non-trivial case, that of two-dimensional, i. e. planar differential systems, cycles may be isolated (limit cycles) or belong to continua of cycles, usually called period annuli. If the inner boundary of a period annulus consists of a single point, then such a point is said to be a center. If (\ref{sysn}) has a period annulus $A$, one can define the period function, by assigning to every $z \in A$ the minimal positive period of  $\gamma_z(t)$. The period function $T(z)$ is a differentiable first integral of  (\ref{sysn}) defined on $A$. Its properties are strictly related to the existence an multiplicity of solutions to some Dirichlet,  Neumann and mixed BVP's \cite{Ca,MMP}, in particular when the system is equivalent to a second order ODE. The period function has been the object of several papers in the last decades \cite{BBLT, CLH, CS, GV, LZ, L, LH}.

In this paper we prove that the study of the cycles' periods in system (\ref{sysn})  can be reduced to that of the cycles' periods in another system having the same orbits, but not the same solutions,
\begin{equation} \label{sysna}
\dot z = \alpha(z) V(z), \qquad \alpha \in \C^0(\Omega,\R).
\end{equation}  
In some cases the second system has a simpler form, allowing for an easier study. In particular for a class of polynomial systems the second system has lower degree than the first one, as shown in next sections' corollaries.

In order to find such simpler systems, we assume (\ref{sysn}) to satisfy a suitable symmetry condition. We consider two kind of symmetries. We assume the existence of an involution $\sigma\in C^0(\Omega,\Omega)$, such that the local flow satisfies either 
$$
\sigma(\phi(t,z)) = \phi(-t,\sigma(z))
$$
($\sigma$-reversibility) or
$$
\sigma(\phi(t,z)) = \phi(t,\sigma(z)).
$$
($\sigma$-symmetry). In other words, we  require the flow to be either anti-commutative or commutative with the symmetry $\sigma$. 

If $\sigma \in C^1(\Omega,\Omega)$, then the $\sigma$-reversibility of (\ref{sysn}) can be checked by verifying the relationship
$$
 V(\sigma(z)) = - J_\sigma(z) \cdot V(z) ,
$$
where $J_\sigma(z)$ is the Jacobian matrix of $\sigma$ at $z$. On the other hand, if $\sigma \in C^1(\Omega,\Omega)$, then the  $\sigma$-symmetry is equivalent to
$$
  V(\sigma(z)) = J_\sigma(z) \cdot V(z).
$$
In the study of planar systems, reversibility has been often considered in order to prove the existence of a center \cite{Ch,GM}. In \cite{C}, sufficient conditions for reversibility were given. In \cite{BM,GM} the converse problem was considered, i. e. under which conditions a center is reversible in a generalized sense. In \cite{BRT,LZ,L,LH} reversible centers were considered in relation to bifurcation problems. 

\bigskip

We prove that, if $\alpha$ satisfies the condition
\begin{equation}      \label{condalfa}
\alpha(z) + \alpha(\sigma (z) )= 2 \alpha(z)  \alpha(\sigma (z) ),
 \end{equation}
then a $\sigma$-invariant cycle has the same period both in system (\ref{sysn}) and in (\ref{sysna}). Moreover, if  (\ref{sysn}) is planar and has a center $O$ with a $\sigma$-invariant central region, then under condition (\ref{condalfa}) the period functions of (\ref{sysn}) and  (\ref{sysna}) are equal. 

On the other hand, if the system (\ref{sysn}) is symmetric w. r. to a point, as it happens when the vector field has odd components, then it can have limit cycles. Even in such a case, the condition (\ref{condalfa}) implies that a $\sigma$-symmetric cycle has the same period in both systems. 

We prove our main result for $n$-dimensional systems, since the proof is not based on the plane's topological properties. On the other hand, our main interest for applications lies in the study of $T(z)$, hence we apply such results to planar systems.

\section{$\sigma$-reversible systems}

\begin{definition} \label{defrev}The local flow $\phi(t,z)$ is said to be {\it $\sigma$-reversible} in $\Omega$ if there exists a map $\sigma\in C^0(\Omega,\Omega)$ such that
\begin{itemize}
\item[(i)] $\sigma^2 = Id$, i. e. $\sigma(\sigma(z)) = z$ for all $z\in \Omega$;
\item[(ii)] $\sigma(\phi(t,z)) = \phi(-t,\sigma(z))$, for all $z\in \Omega$ such that both sides of such equality are defined.
\end{itemize}
A differential system (a vector field) is said to be {\it $\sigma$-reversible} if its local flow is e {\it $\sigma$-reversible}.
\end{definition}

We say that a point $z$ is a $\sigma$-fixed point if $\sigma(z) = z$. Two points $z^+$, $z^{++}$ are said to be $\sigma$-related if 
$\sigma(z^+) = z^{++}$. Such a relationship is symmetric, since $\sigma$ is an involution. 

The above definition extends the usual definition of reversible system, in which $\sigma$ is assumed to be the mirror symmetry w. r. to a line. Such a system orbit set is itself symmetric with respect to a line. In particular, if an orbit intersects the symmetry line at two distinct points, then it is a cycle. Such systems have been studied by several authors, in particular in relation to the existence of centers. A special case is that of planar vector fields whose first component is $x$-even, and the second component is $x$-odd,
\begin{equation} \label{xrev}
V_1(-x,y) = V_1(x,y) , \qquad V_2(-x,y) = -V_2(x,y).
 \end{equation}
In this case the orbits are symmetric with respect to the $y$-axis. If the origin $O$ is a critical point of rotation type, then it is a center. A planar differential system may be reversible with respect to several lines. For instance, the system

$$ 
\left\{
\matrix{ \dot x &= y  \hfill  \cr \dot y  &= - x   ,  }
\right.
$$
is reversible w. r. to every line through $O$, and the system
$$ 
\left\{
\matrix{ \dot x &= y(1+x^2)  \hfill  \cr \dot y  &= - x(1+y^2)   ,  }
\right.
$$
is reversible w. r. both to the $y$-axis and to the $x$-axis.

In \cite{C} conditions to determine all  the possible reversibility lines of non-hamiltonian systems are studied. 

\begin{definition}Let $U$ be an open subset of $\Omega$, $\sigma  \in C^0(U , U)$, $\alpha \in C^0(U,\R)$. Then $\alpha$ is said to be $\sigma $-compatible on $U$ if, for all $z \in U$, the equality (\ref{condalfa}) holds.
 \end{definition}

 If $\alpha(z) \neq 0$ is $\sigma$-compatible, setting $u= \alpha(z)$, $v = \alpha(\sigma(z))$, we see that the point$(u,v)=(\alpha(z),\alpha(\sigma(z)))$ lies on the  branch of the equilateral hyperbola
$$
2uv -u -v = 0
$$
contained in the first orthant. If $z^*$ is a $\sigma$-fixed point, one has 
$$
2\alpha(z^*) ^2 - 2\alpha(z^*)  = 0,
$$
hence $\alpha(z^*) = 1$. 

\begin{definition}
We say that a function $\delta \in C^0(\Omega,\Omega)$ is $\sigma$-odd if
$$
\delta(\sigma(z)) = - \delta (z).
$$
\end{definition}

We can characterize $\sigma$-compatible functions by means of $\sigma$-odd functions.

\begin{lemma} \label{lemmadelta}The function $\alpha(x,y)$ is $\sigma$-compatible if and only if there exists a $\sigma$-odd function $\delta$ such that 
\begin{equation} \label{ad}
\alpha(z) = \frac {1}{1 + \delta(z)}.
\end{equation}  
\end{lemma}
{\it Proof.}   One may divide the equality
$$
\alpha(z) + \alpha (\sigma(z)) = 2 \alpha(z)  \alpha (\sigma(z)) 
$$
by the product $\alpha(z)  \alpha (\sigma(z))$, obtaining
$$
\frac{1}{\alpha(z)} + \frac{1}{\alpha(\sigma(z))} = 2 
$$
\begin{equation} \label{a}
\frac{1}{\alpha(z)} - 1 =  1 -  \frac{1}{\alpha(\sigma(z))} = -\left( \frac{1}{\alpha(\sigma(z))} - 1 \right),
 \end{equation}
hence the function
$$
\delta(z) = \frac{1}{\alpha(z)} - 1 
$$
is $\sigma$-odd. 

Such steps may be read backwards, proving that if $\delta$ is $\sigma$-odd, then (\ref{ad}) implies (\ref{condalfa}).
\hfill  $\clubsuit$ \\

The sign condition $\alpha(z) > 0$ is equivalent to $\delta(z) > -1$. 

\bigskip

Let $\alpha \in C^0(\Omega,\R)$ be a  function such that $\alpha(z) \neq 0$ for all $z \neq O$. The system (\ref{sysna})
has the same orbits as (\ref{sysn}). In general, such systems do not have equal period functions. 
We denote by $\phi_\alpha(t_\alpha,z)$ the local flow defined by the solutions of (\ref{sysna}).

\begin{lemma} \label{lemmarev}
Let (\ref{sysn})  be $\sigma$-reversible and $\gamma$ be a $T$-periodic cycle of (\ref{sysn}).  If $\gamma$ contains two distinct related points, then $\gamma$ contains two distinct $\sigma$-fixed points. Moreover, if $z^*$ is a $\sigma$-fixed point, then also $\phi\left( \frac T2,z^* \right)$ is a $\sigma$-fixed point.
\end{lemma}
{\it Proof.}   
Let $T$ be the (minimal positive) period of $\gamma$, $z^+$, $z^{++}$ the distinct $\sigma$-related points on $\gamma$.
$\gamma$ is a cycle, hence there exists $t \in (0,T)$ such that $\phi(t,z^+) = \sigma(z^+) =z^{++} $. Let us set $z^* = \phi\left( \frac t2,z^+ \right)$. Then
$$
\sigma(z^*) = \sigma \left( \phi\left( \frac t2,z^+ \right) \right) = 
  \phi\left( -\frac t2, \sigma(z^+) \right) =   \phi\left( -\frac t2, z^{++}  \right) =
  $$
  $$=    \phi\left( -\frac t2, \phi(t,z^+)  \right) = \phi \left (\frac t2,z^+ \right) = z^*.
$$
 
In order to prove the existence of a second $\sigma$-fixed point, let $z^*$ a  $\sigma$-fixed point : $\sigma(z^*)=z^*$. Let us set $z^{**} = \phi \left(\frac T2,z^*\right)$. Then one has, by point (ii) of  definition  \ref{defrev} and by the $T$-periodicity of $\gamma$,
$$
\sigma(z^{**}) = \sigma\left( \phi\left(\frac T2,z^*\right)\right) = \phi\left(-\frac T2,\sigma( z^*)  \right) 
=  \phi\left(-\frac T2, z^*  \right) =  \phi\left(\frac T2, z^*  \right) = z^{**}.
$$
Since $z^* \neq z^{**}$, $\gamma$ contains two distinct $\sigma$-fixed points.
\hfill  $\clubsuit$ \\

It is easy to prove that a cycle $\gamma$ contains two distinct related points if and only if $\sigma(\gamma) = \gamma$.
\begin{theorem} \label{teorev}
Let (\ref{sysn})  be $\sigma$-reversible and $\gamma$ be a cycle of (\ref{sysn}) containing two $\sigma$-related points.  If $\alpha \in C^0(\Omega,\R)$ is  $\sigma$-compatible on a neighbourhood $U$ of $\gamma$ and $\alpha (z) \neq 0$ on $\gamma$, then the period of $\gamma$ is the same both for (\ref{sysn}) and (\ref{sysna}). 
\end{theorem}
{\it Proof.}   
Without loss of generality we may assume $\alpha(z) > 0$ for $z \neq O$. 
Let us denote by $T$ $\gamma$'s period with respect to (\ref{sysn}), and $T_\alpha$ $\gamma$'s period with respect to (\ref{sysna}). One has
$$
T = \int_0^T dt , \qquad T_\alpha = \int_0^{T_\alpha} dt_\alpha,
 $$
where $t$ and $t_\alpha$are the time-variables of the two systems and the integration is performed along the cycles $\phi(t,z)$ and $\phi_\alpha(t_\alpha,z)$. Since $\alpha >0$, along $\gamma$ we may write $t$ as a function of $t_\alpha$ and viceversa. Such functions is differentiable and one has
$$
\frac{d t}{d t_\alpha} = \alpha, \qquad \frac{d t_\alpha} {d t}= \frac 1\alpha.
$$
Moreover, changing the integration variable in the second integral above,
$$
 T_\alpha = \int_0^{T_\alpha} dt_\alpha  = 
 \int_0^{T}   \frac 1\alpha dt =  \int_0^{T} \left( 1 - 1 + \frac 1\alpha  \right)dt = T+ \int_0^{T} \left(  \frac 1\alpha -1  \right)dt .
$$
Hence the two periods coincide if and only if 
$$
 \int_0^{T} \left(  \frac 1\alpha -1  \right)dt = 0.
 $$
By lemma \ref{lemmarev} one $\gamma$ contains a $\sigma$-fixed point $z^*$. The above integral can be computed along the solution $\phi(t,z^*)$. We prove the above equality by showing that
$$
 \int_0^{\frac T2} \left(  \frac 1{\alpha\left( \phi(t,z^*) \right) }-1  \right)dt =
  - \int_{\frac T2}^T  \left( \frac 1{\alpha\left( \phi(t,z^*) \right) }-1 \right) dt .
 $$
Writing $z^{**} = \phi\left( \frac T2, z^* \right)$, by lemma \ref{lemmarev} $\sigma(z^{**}) = z^{**}$. Then one has 
$$
- \int_{\frac T2}^T  \left( \frac 1{\alpha\left( \phi(t,z^*) \right) }-1 \right) dt =
-\int_{\frac T2}^T  \left( \frac 1{\alpha\left( \phi(t - \frac T2,z^{**}) \right) }-1 \right) dt = 
$$
$$
-\int_0^{\frac T2}  \left( \frac 1{\alpha\left( \phi(t ,z^{**}) \right) }-1 \right) dt = 
-\int_0^{\frac T2}  \left( \frac 1{\alpha\left( \phi(t ,\sigma(z^{**})) \right) }-1 \right) dt = 
$$
$$
-\int_0^{\frac T2}  \left( \frac 1{\alpha(\sigma \left( \phi(-t ,z^{**}) )\right) }-1 \right) dt =
-\int_0^{\frac T2}  \left(1 - \frac 1{\alpha \left( \phi(-t ,z^{**})\right) } \right) dt .
$$
by the equality (\ref{a}). Now, changing the integration variable, $\tau = -t$, one has 
$$
-\int_0^{\frac T2}  \left(1 - \frac 1{\alpha \left( \phi(-t ,z^{**})\right) } \right) dt  =
- \int_{- \frac T2}^0  \left( 1  -\frac 1{\alpha \left( \phi(\tau ,z^{**})\right)  }  \right) d \tau =
$$
$$
- \int_{- \frac T2}^0  \left( 1  -\frac 1{\alpha \left( \phi(\tau + \frac T2,z^*)\right)  }  \right) d \tau =
- \int_0^{\frac T2}  \left( 1  -\frac 1{\alpha \left( \phi(s,z^*)\right)  }  \right) d s =
$$
$$
= \int_0^{\frac T2}  \left( \frac 1{\alpha \left( \phi(s,z^*)\right)  } -1  \right) d s ,
$$
as required.
\hfill  $\clubsuit$ \\

An application to a specific class of systems is given by next corollary.
\begin{corollary} \label{correv} Let  the system
$$
\dot z = V(z)(\delta(z)^2 - 1).
$$
be $\sigma$-reversible, and $\delta(z)$ be $\sigma$-odd, with $\delta(z) > -1$ on $\Omega$. If $\gamma$ is a cycle containing two $\sigma$-related points, then its period in the above system and in the system
$$
\dot z = V(z)(\delta(z) - 1)
$$
is the same.
\end{corollary}
{\it Proof.}   The former system satisfies the hypotheses of theorem \ref{teorev}. The latter system is obtained from it  multiplying by $\displaystyle{ \alpha(z) = \frac 1{1 + \delta(z)}   }$. By  theorem \ref{teorev}, they have the same period function.
\hfill  $\clubsuit$ \\

We restrict to the simplest reversibility condition, i. e. that one obtained by choosing $\sigma(x,y) = (-x,y)$. Such reversibility can be verified by checking the conditions (\ref{xrev}). Under such a choice, a function $\alpha$ is $\sigma$-compatible if and only if
$$
\alpha(x,y) + \alpha(-x,y)= 2 \alpha(x,y)  \alpha(-x,y).
$$
Moreover, a function $\delta(x,y)$ is $\sigma$-odd if and only if  $\delta(-x,y) = - \delta(x,y)$, i. e. if it is $x$-odd.
We can construct a wide class of systems with equal period functions perturbing hamiltonian systems with separable variables. If $f(y)$, $g(x)$, are odd, with $yf(y) > 0$ for $y \neq 0$, $xg(x) > 0$ for $x \neq 0$, $\delta(x,y)$ is $x$-odd and $F(x,y)$ is $x$-even, non-zero in a punctured neighbourhood of $O$, then the systems
$$ 
\left\{
\matrix{ \dot x &= f(y)\Big(1 - \delta(x,y)^2\Big)F(x,y)   \hfill  \cr \dot y  &= - g(x)\Big(1 - \delta(x,y)^2\Big) F(x,y),  }
\right.
$$
and
$$
 \left\{
\matrix{ \dot x &= f(y)\Big(1 - \delta(x,y)\Big) F(x,y)  \hfill  \cr \dot y  &= - g(x)\Big(1 - \delta(x,y)\Big)F(x,y) ,  }
\right.
$$
have both a center at $O$ and equal period functions. We get a very simple example of such a system by taking $f(y) = y$, $g(x) = x$, $\delta(x,y) = x$, $F(x,y) = 1$. The systems
$$
 \left\{
\matrix{ \dot x &= y(1 - x^2)    \hfill  \cr \dot y  &= - x(1 - x^2 ) ,  }
\right.
$$
and
$$ 
\left\{
\matrix{ \dot x &= y (1 - x  )    \hfill  \cr \dot y  &= - x (1 - x  )  ,  }
\right.
$$
have equal period functions.

We can consider also systems exhibiting two types of reversibility. For instance,
the system
$$
 \left\{
\matrix{ \dot x &= y(1 - x^2) (1-y^4)   \hfill  \cr \dot y  &= - x(1 - x^2 )(1-y^4) ,  }
\right.
$$
is reversible both w. r. to the $x$-axis and w. r. to the $y$-axis, the symmetries being given by $\sigma^y(-x,y) = (-x,y)$, $\sigma^x (x,y) = (x,-y)$. We can divide the vector field either by $\displaystyle{ 1 + \delta^y(x,y) = 1 + x  }$, or by $\displaystyle{ 1 + \delta^x(x,y) =   1 - y  }$, obtaining two systems whose period functions are equal and coincide with that of the original system,
 $$
 \left\{
\matrix{ \dot x &= y(1 - x) (1-y^4)   \hfill  \cr \dot y  &= - x(1 - x )(1 - y^4) ,  }
\right.
$$
$$
 \left\{
\matrix{ \dot x &= y(1 - x^2) (1+y+y^2+y^3)   \hfill  \cr \dot y  &= - x(1 - x^2 )(1+y+y^2+y^3) .  }
\right.
$$

\section{$\sigma$-symmetric systems}

\begin{definition} \label{defsigma}The local flow $\phi(t,z)$ is said to be {\it $\sigma$-symmetric} if there exists a map $\sigma\in C^0(\Omega,\Omega)$ such that
\begin{itemize}
\item[i)] $\sigma^2 = Id$, i. e. $\sigma(\sigma(z)) = z$ for all $z\in \Omega$;
\item[ii)] $\sigma(\phi(t,z)) = \phi(t,\sigma(z))$, for all $(t,z) \in \R \times \Omega$ such that both sides of such equality are defined.
\end{itemize}
A differential system is said to be {\it $\sigma$-symmetric} if its local flow is  {\it $\sigma$-symmetric}.
  \end{definition}

The above definition extends the usual definition of symmetry with respect to a point. Such systems do not necessarily have centers, even if an isolated critical point is of rotation type. They may have  attractive or repulsive critical points and cycles, i. e. limit cycles, as the system
$$
\dot x = y , \qquad \dot y =  -x -y( x^2 - 1),
$$
equivalent to Van der Pol equation and symmetric with respect to the origin.

\begin{lemma} \label{lemmarsimm}
Let (\ref{sysn})  be $\sigma$-symmetric and $\gamma$ be a $T$-periodic cycle of (\ref{sysn}) containing two $\sigma$-related points, $z^{+}$ and $z^{++}$. Then for all $z \in \gamma$
$$
\sigma(z) = \phi \left(\frac T2, z \right) .
$$
\end{lemma}
{\it Proof.}   
Let $T$ be the minimal positive period of $\gamma$.
There exists $t\in (0,T)$ such that $\phi(t,z^{+}) = \sigma(z^{+}) = z^{++}$.
Possibly exchanging $z^{+}$ and $z^{++}$, we may assume $t\in \left(0,\frac T2 \right]$.
One has
$$
z^{+} = \sigma^2(z^{+}) = \sigma(\sigma(z^{+}) )= \sigma(\phi(t,z^{+})) = \phi(t,\sigma(z^{+})) = \phi(t,\phi(t,z^{+})) = \phi(2t,z^{+}). 
$$
Since $0 < 2t \leq T$, one has $T = 2t$, i. e. $\sigma(z^{+}) = \phi \left(\frac T2, z^{+} \right) $. Then, for an arbitrary $z \in \gamma$ there exists $\tau \in (0,T)$ such that $z=\phi(\tau,z^{+})$. Hence
$$
\sigma (z) = \sigma \left(\phi \left(\tau, z^{+} \right)\right) = \phi\left( \tau, \sigma(z^{+}) \right) = 
 \phi\left( \tau,  \phi \left(\frac T2, z^{+} \right)  \right) =
$$
$$ =  \phi \left(\tau + \frac T2, z^{+} \right)
=  \phi \left( \frac T2, \phi( \tau,z^{+}) \right) =   \phi \left(\frac T2, z \right) .
$$
\hfill  $\clubsuit$ \\

The above lemma shows that a $\sigma$-symmetric flow does not have non-stationary $\sigma$-fixed points.

Next theorem is analogous to theorem \ref{teorev}. Its proof is somewhat easier, since it does not require to select two special points of a cycle $\gamma$, fixed w. r. to $\sigma$. Actually, if the system is $\sigma$-symmetric, a cycle contain no $\sigma$-fixed points. As in the case of $\sigma$-reversible systems, a cycle $\gamma$ contains two distinct related points if and only if $\sigma(\gamma) = \gamma$.

\begin{theorem} \label{teosimm}
Let (\ref{sysn})  be $\sigma$-symmetric and $\gamma$ be a cycle of (\ref{sysn}) containing two $\sigma$-related points.  If $\alpha \in C^0(\Omega,\R)$ is  $\sigma$-compatible on a neighbourhood $U$ of $\gamma$ and $\alpha (z) \neq 0$ on $\gamma$, then the period of $\gamma$ is the same both for (\ref{sysn}) and (\ref{sysna}). 
\end{theorem}
{\it Proof.}   Without loss of generality we may assume $\alpha(z) > 0$ for $z \neq O$. 
The first part of the proof is identical to that of  theorem \ref{teorev}. Let $z^*$ be an arbitrary point of $\gamma$. As in theorem \ref{teorev}, it is sufficient to prove that 
$$
 \int_0^{T} \left(  \frac 1{\alpha(\phi(t,z^*)) }-1  \right)dt  =  0.
 $$
One has
 $$
 \int_{\frac T2}^T \left(  \frac 1{\alpha(\phi(t,z^*)) }-1  \right)dt =
  \int_{\frac T2}^T \left(1 -  \frac 1{\alpha (\sigma(\phi(t,z^*))) }  \right)dt =
 $$
$$
  \int_{\frac T2}^T \left(1 -  \frac 1{\alpha (\phi(t,\sigma(z^*))) }  \right)dt =
  \int_{\frac T2}^T \left(1 -  \frac 1{\alpha (\phi(t+\frac T2,z^*))) }  \right)dt =
 $$
$$
  \int_T ^{\frac {3T}2}\left(1 -  \frac 1{\alpha (\phi(u,z^*))) }  \right)du =
  \int_0 ^{\frac {T}2}\left(1 -  \frac 1{\alpha (\phi(u,z^*))) }  \right)du ,
 $$
since $\gamma$ is $T$-periodic. Then
$$
 \int_0^{T} \left(  \frac 1{\alpha(\phi(t,z^*)) }-1  \right)dt  =  
 $$  
 $$
  \int_0^{\frac T2} \left(  \frac 1{\alpha(\phi(t,z^*)) }-1  \right)dt  +  \int_{\frac T2}^T \left(  \frac 1{\alpha(\phi(t,z^*)) }-1  \right)dt =
 $$
$$
  \int_0^{\frac T2} \left(  \frac 1{\alpha(\phi(t,z^*)) }-1  \right)dt  +  \int_0 ^{\frac {T}2}\left(1 -  \frac 1{\alpha (\phi(t,z^*))) }  \right)dt = 0.
 $$
\hfill  $\clubsuit$ \\

The proof of next corollary is similar to that of corollary \ref{correv}.

\begin{corollary} \label{corrsimm} Let  the system
$$
\dot z = V(z)(\delta(z)^2 - 1).
$$
be $\sigma$-symmetric, and $\delta(z)$ be $\sigma$-odd, with $\delta(z) > -1$ on $\Omega$. If $\gamma$ is a cycle containing two $\sigma$-related points, then its period in the above system and in the system
$$
\dot z = V(z)(\delta(z) - 1)
$$
is the same.
\end{corollary} 

We restrict to the simplest  symmetry condition, i. e. that one obtained by choosing $\sigma(z) = -z$. Such a condition can be easily verified by checking the condition 
\begin{equation} \label{simm}
V(-z) = -V(z).
 \end{equation}
In this case the orbit set is symmetric with respect to a point, the center of symmetry. Rotation points are not necessarily centers, limit cycles may exist. 

\bigskip

Choosing $\sigma(z) = -z$, a function $\delta(x,y)$ is $\sigma$-odd if and only if  $\delta(-z) = - \delta(z)$.

\bigskip

Working similarly to the previous section, we can construct a wide class of systems with equal period functions perturbing hamiltonian systems with separable variables. If $f(y)$, $g(x)$, are odd, with $yf(y) > 0$ for $y \neq 0$, $xg(x) > 0$ for $x \neq 0$, $\delta(-x,-y) = - \delta(x,y)$, $F(-x,-y)= F(x,y)$, non-vanishing in a punctured neighbourhood of $O$, then the systems
$$ 
\left\{
\matrix{ \dot x &= f(y)\Big(1 - \delta(x,y)^2\Big)F(x,y)   \hfill  \cr \dot y  &= - g(x)\Big(1 - \delta(x,y)^2\Big) F(x,y),  }
\right.
$$
and
$$
 \left\{
\matrix{ \dot x &= f(y)\Big(1 - \delta(x,y)\Big) F(x,y)  \hfill  \cr \dot y  &= - g(x)\Big(1 - \delta(x,y)\Big)F(x,y) ,  }
\right.
$$
have both a center at $O$ and equal period functions. For instance, we may take $f(y) = y$, $g(x) = x$, $\delta(x,y) = x+y$, $F(x,y) = 1$. The systems
$$
 \left\{
\matrix{ \dot x &= y(1 - (x+y)^2)    \hfill  \cr \dot y  &= - x(1 -  (x+y)^2 ) ,  }
\right.
$$
and
$$ 
\left\{
\matrix{ \dot x &= y (1 - x  -y)    \hfill  \cr \dot y  &= - x (1 - x -y )  ,  }
\right.
$$
have equal period functions.

\bigskip

Similarly, we may consider a system with a limit cycle, equivalent to a Li\'enard equation:
$$ 
\left\{
\matrix{ \dot x &= y  \hfill  \cr \dot y  &= - x - yf(x)   ,  }
\right.
$$
We assume $f(x)$ to be even, $f(-x) = f(x)$, so that the system is symmetric w. r. to the origin.
If $\delta(x)$ is odd, then the two systems
\begin{equation}  \label{equa1}
\left\{
\matrix{ \dot x &= y (1 - \delta(x)^2)  \hfill  \cr \dot y  &= ( - x - yf(x))   (1 - \delta(x)^2)  ,  }
\right.
\end{equation}
and
\begin{equation}  \label{equa2}
\left\{
\matrix{ \dot x &= y (1 - \delta(x))  \hfill  \cr \dot y  &= (- x - yf(x) ) (1 - \delta(x))  ,  }
\right.
\end{equation}
in the set $\{ (x,y) : \delta (x) > -1   \}$ have the same cycles with the same periods. Even if both systems do not appear in the traditional form equivalent to a second order ODE,  both are equivalent to a second order ODE. The system (\ref{equa1}) is equivalent to
$$
\ddot x + \dot x (1 - \delta(x) )f(x) + \dot x^2 \frac{\delta'(x)}{1 - \delta(x)} + x (1 - \delta(x))^2 = 0.
$$
The system (\ref{equa2}) is equivalent to
$$
\ddot x + \dot x (1 - \delta(x)^2 )f(x) + \dot x^2 \frac{2 \delta(x)\delta'(x)}{1 - \delta(x)^2} + x( (1 - \delta(x))^2 )^2= 0.
$$
This implies that one can establish a period-preserving one-to-one correspondence between periodic solutions of both systems.

\enddocument